\documentclass[12pt,
reqno]{amsart}

\usepackage[cp1251]{inputenc}
\usepackage{amsmath,amsxtra,amssymb,eufrak}
\usepackage[all]{xy}
\usepackage{mathrsfs}
\usepackage[english]{babel}

\theoremstyle{plain}
\newtheorem{thm}{Theorem}[section]
\newtheorem{lm}[thm]{Lemma}
\newtheorem{co}[thm]{Corollary}
\newtheorem{pr}[thm]{Proposition}

\theoremstyle{definition}

\newtheorem{rem}[thm]{Remark}
\newtheorem{que}[thm]{Question}

\numberwithin{equation}{section}

\DeclareMathOperator{\Ker}{Ker}
\newcommand*{\ptens}[1]{\mathop{\widehat\otimes}_{#1}}
\renewcommand{\le}{\leqslant}
\renewcommand{\ge}{\geqslant}
\title{Homological dimensions of algebras of analytic functionals and their completions}

\author{O. Yu. Aristov}
\address{Institute for Advanced Study in Mathematics of Harbin Institute of Technology, Harbin 150001, China;
\newline\indent
Suzhou Research Institute of Harbin Institute of Technology, Suzhou 215104, China}
\email{aristovoyu@inbox.ru}

\begin{document}

\begin{abstract}
We show that the main homological dimensions of the algebra of analytic functionals on a connected complex Lie group, as well as some of its completions, coincide with the dimension of the simply connected solvable factor in the canonical decomposition of the linearization of this group. Thus, the possible nontriviality of a linearly complex reductive factor does not affect the homological properties of the algebras under consideration.
\end{abstract}

\thanks{Translated by DeepSeek and the author}

 \maketitle
 \markright{Homological dimensions}

\section{Introduction}

In \cite{AHHFG, Ar_smash,Ar_aahe}, the author has been studied properties of a classical object of noncommutative harmonic analysis --- the algebra of analytic functionals ${{\mathscr A}}(G)$ on a connected complex Lie group~$G$ --- as well as some completions of this algebra. In this article, the results obtained earlier are used to compute homological dimensions (projective, as well as weak, i.e., flat).

The study of homological dimensions of Banach algebras within the framework of topological homology was initiated by A.\,Ya.~Helemskii in the 1970s. Computing and estimating of these invariants in the Banach case turned out to be exceptionally difficult, but some progress has been made; see the book \cite{X1}, the later survey \cite{He00}, as well as papers \cite{Ar00,Se04,Se07,Se18} not included in the survey. Somewhat later, J.~Taylor, motivated by questions of spectral theory, investigated homological problems for more general topological algebras \cite{T0,T2}; see also \cite{He81}. It turned out that nuclear algebras are better suited to topological homology than Banach algebras, which is due to the properties of the projective tensor product. (The results proved here confirm this observation --- all algebras discussed below are nuclear.) Homological dimensions of nuclear function algebras were computed in the $C^\infty$-case in \cite{Og86,OH84}, and in the holomorphic case in \cite{Pi14A} (for commutative algebras) and \cite{Do03,Pi09} (for noncommutative algebras). Among more recent works, we would like to mention \cite{Ko21,Ko23}. The articles \cite{Pi08,Pi12,Se96,Se25+}, where the weak dimensions are discussed, are also of interest; see in addition \cite{Pi02,Pi10,Pi11A,PS07}. For convolution algebras, these invariants have been studied mainly in the Banach case; see, e.g., \cite{X1,Se18}. An exception is Taylor's result on algebras of distributions on compact real Lie groups \cite{T2} mentioned below. But the homological dimensions of algebras of analytic functionals have not  been previously explored, let alone their completions, which represent an entirely new area of research.

A classical  theorem of Hochschild states that a connected linear complex Lie group admits a decomposition into a semidirect product $B\rtimes L$, where $B$ is simply connected and solvable, and $L$ is linearly complex reductive \cite{HO65}; see also \cite[p.\,601, Theorem 16.3.7]{HiNe}. In particular, if~$\Lambda$ is the linearizer of a connected Lie group~$G$, then the decomposition is possible for $G/\Lambda$. Here we show that \emph{the main homological dimensions of the algebra ${{\mathscr A}}(G)$ coincide with the dimension of~$B$} (this is part of the statement of Theorem~\ref{homdimth}).

The remaining statements in the above mentioned theorem concern some completions of ${{\mathscr A}}(G)$, which are Fr\'echet--Arens--Michael algebras. The Arens--Michael envelope $\widehat{{\mathscr A}}(G)$  and $\widehat{{\mathscr A}}(G)^{\mathrm{PI}}$ --- the envelope with respect to the class of Banach PI-algebras --- are of particular interest to us  (the first envelope has been known for a long time, and the second has been introduced in~\cite{Ar_smash}). The algebra $\widehat{{\mathscr A}}(G)$ can be defined as the completion of ${{\mathscr A}}(G)$ with respect to the family of all continuous submultiplicative seminorms, and $\widehat{{\mathscr A}}(G)^{\mathrm{PI}}$ as the completion with respect to the family of continuous submultiplicative seminorms with the additional restriction: the completion with respect to each of them is an (automatically Banach) algebra satisfying a polynomial identity.

Also,  a unified approach to both algebras has been proposed in the aforementioned article \cite{Ar_smash} --- one can considers them as extreme cases of a certain family. The members of this family are denoted by ${\mathscr A}_{\omega_{max}^\infty}(G/\Lambda)$, where $\omega_{max}$ is a submultiplicative weight depending on the choice of a nilpotent subgroup of $G/\Lambda$ intermediate between the exponential and nilpotent radicals. Namely, $\omega_{max}$ is the maximal submultiplicative weight with exponential distortion on this subgroup. In~\cite{Ar_aahe}, it was demonstrated that the map ${{\mathscr A}}(G)\to {\mathscr A}_{\omega_{max}^\infty}(G/\Lambda)$ has an important property --- it is a homological epimorphism; see Theorem~\ref{C4mainexdi}. One of the direct applications of this result is  computing of homological dimensions of completions of the indicated type. This is possible because the property of being a homological epimorphism simplifies the estimation of upper bounds. (The reader interested only in ${{\mathscr A}}(G)$ can easily notice that the arguments relating to it do not include references to homological epimorphisms.) As a special case, we show that \emph{the homological dimensions of $\widehat{{\mathscr A}}(G)$ and $\widehat{{\mathscr A}}(G)^{\mathrm{PI}}$ also coincide with the dimension of~$B$} (see Corollaries~\ref{hdAM} and~\ref{hdPIE}). As for the lower bound, the argument reduces to the question of Lie algebra cohomology, for which the answer was given in~\cite{Ha70}.

Note that  the results obtained imply that the possible non-triviality of the linearly complex reductive factor $L$ in the decomposition $G/\Lambda=B\rtimes L$ does not affect the homological dimensions of ${\mathscr A}(G)$. However, the fact that ${{\mathscr A}}(L)$ is homologically trivial is not sufficient to derive this property --- it is necessary to invoke the statement that ${\mathscr A}(G/\Lambda)$ is \emph{relatively homologically trivial} over ${\mathscr A}(B)$; see details below.

In conclusion,  we also formulate a question about real Lie groups. Denote by ${\mathscr E}'(G)$ the algebra of compactly supported distributions on a connected real Lie group $G$, i.e., the algebra of continuous linear functionals on the space of all smooth functions with convolution as multiplication.
\begin{que}\label{HdEGpr}
What are the homological dimensions of ${\mathcal{E}}'(G)$, $\widehat{\mathcal{E}'}(G)$, and $\widehat{\mathcal{E}'}(G)^{\mathrm{PI}}$?
\end{que}
The only result in this direction known to the author  is the following statement: if $G$ is compact, then $\mathcal{E}'(G)$ is homologically trivial \cite[Proposition 7.3]{T2}.

\section{Preliminary information and auxiliary statements}

We consider $\mathbin{\widehat{\otimes}}$-algebras and $\mathbin{\widehat{\otimes}}$-modules, i.e., complete locally convex algebras and modules with jointly continuous multiplication (over the ground field $\mathbb{C}$). All algebras and modules, as well as their (homo)morphisms, are assumed to be unital. Furthermore, we need $\mathbin{\widehat{\otimes}}$-coalgebras, $\mathbin{\widehat{\otimes}}$-bialgebras, and $\mathbin{\widehat{\otimes}}$-Hopf algebras, i.e., coalgebras, bialgebras, and Hopf algebras in the symmetric monoidal category of complete locally convex spaces over~$\mathbb{C}$ endowed  with the bifunctor $(-)\mathbin{\widehat{\otimes}} (-)$ of the complete projective tensor product; see, e.g.,~\cite{Ar_smash}.

\subsection*{Algebras of analytic functionals and their completions}
For a complex Lie group $G$, we consider its algebra of analytic functionals ${\mathscr A}(G)$, i.e., the space of continuous functionals on the algebra of holomorphic functions on $G$ endowed with convolution. (In the generality we need, they first appeared in~\cite{Li70}.)  The completion corresponding a submultiplicative weight $\omega$ on $G$ is a Banach algebra \cite{Ak08}. Moreover, using all powers of $\omega$, we obtain the Fr\'echet algebra ${\mathscr A}_{\omega^\infty}(G)$, which is an Arens--Michael algebra, i.e., a projective limit of Banach algebras. For details  see~\cite{Ar_smash}.

\subsection*{Analytic smash products}
In what follows, it is essential that algebras of analytic functionals are also $\mathbin{\widehat{\otimes}}$-Hopf algebras~\cite{Ar_smash}. This allows us to turn to a analytic version of the smash product construction, which was considered by Pirkovskii in~\cite{Pi4}. It was shown in~\cite[\S\,2]{Ar_smash} that semidirect products of Lie groups correspond to analytic smash products of algebras of analytic functionals.

Let us give the necessary definitions. Suppose that $H$ is a $\mathbin{\widehat{\otimes}}$-Hopf algebra and $A$ is a $\mathbin{\widehat{\otimes}}$-algebra endowed  with a structure of a left $H$-$\mathbin{\widehat{\otimes}}$-module. In this text, we freely use the generalized Sweedler notation for $\mathbin{\widehat{\otimes}}$-Hopf algebras. The possibility of generalizing Sweedler notation was first pointed out by Akbarov in \cite{Ak08}. A more detailed discussion can be found in~\cite[\S\,1]{Ar_smash}. Recall that~$A$ is called a (left) \emph{$H$-$\mathbin{\widehat{\otimes}}$-module algebra} if the conditions
\begin{equation}\label{Hmodalgcond}
h\cdot(ab)=\sum (h_{(1)}\cdot a)(h_{(2)}\cdot b)\quad\text{and}\quad h\cdot 1=\varepsilon(h)1\qquad(h\in H,\,a,b\in A)
\end{equation}
hold. (As usual, $\varepsilon$ denotes the counit.)
Then  the formula
\begin{equation}\label{exmusmpr}
(a\otimes h)(a'\otimes h')=\sum a(h_{(1)}\cdot a')\otimes h_{(2)}h'\qquad(h,h'\in H;\,a,a'\in A)
\end{equation}
defines a jointly continuous associative multiplication on the complete locally convex space $A\mathbin{\widehat{\otimes}} H$. This space endowed  with this multiplication is denoted by $A\mathop{\widehat{\#}} H$ and is called an \emph{analytic smash product}. (Elementary tensors are denoted by $a\mathop{\#} h$.)

Let now $C$ be a $\mathbin{\widehat{\otimes}}$-coalgebra endowed  with a structure of a left $H$-$\mathbin{\widehat{\otimes}}$-module. If
\begin{equation}\label{coumo}
\sum (h\cdot c)_{(1)}\otimes (h\cdot c)_{(2)}=\sum (h_{(1)}\cdot c_{(1)})\otimes(h_{(2)}\cdot c_{(2)})\quad\text{and}\quad \varepsilon_C(h\cdot c)= \varepsilon_H(h)\varepsilon_C(c)
\end{equation}
for all $h\in H$ and $c\in C$, it is called a (left) \emph{$H$-$\mathbin{\widehat{\otimes}}$-module coalgebra}; cf. \cite[Definition 2.1(c)]{Mo77}.
In the case when $A$ is a $\mathbin{\widehat{\otimes}}$-bialgebra that is simultaneously an
$H$-$\mathbin{\widehat{\otimes}}$-module algebra and an $H$-$\mathbin{\widehat{\otimes}}$-module coalgebra,  it
is called a (left) \emph{$H$-$\mathbin{\widehat{\otimes}}$-module bialgebra}; cf. \cite[Definition 2.1(e)]{Mo77}.  For details see~\cite[\S\,1]{Ar_smash}.

\subsection*{Homological dimensions}

Recall basic definitions from the homological theory of locally convex algebras; see, e.g., \cite{X1} and \cite{Pir_qfree}.

Consider an $A$-$\mathbin{\widehat{\otimes}}$-module $P$ (left, right, or bimodule) over a $\mathbin{\widehat{\otimes}}$-algebra $A$. It is said to be \emph{projective} if for any admissible epimorphism of $A$-$\mathbin{\widehat{\otimes}}$-modules (i.e., having a right inverse continuous linear map) with codomain $P$, there exists a right inverse morphism of $A$-$\mathbin{\widehat{\otimes}}$-modules. A special case of a projective module is a \emph{free} module, i.e., one isomorphic to a module of the form $A\mathbin{\widehat{\otimes}} E$ for some complete locally convex space~$E$ (for left modules). A chain complex
$$
\cdots\leftarrow M_n \leftarrow M_{n+1} \leftarrow \cdots
$$
of $A$-$\mathbin{\widehat{\otimes}}$-modules is said to be \emph{admissible} if it is contractible in the category of topological vector spaces, i.e., if there exists a contracting homotopy consisting of continuous linear maps.

If $M$ and $N$ are a right and a left $A$-$\mathbin{\widehat{\otimes}}$-module, respectively, then their \emph{$A$-module tensor product} $M\ptens{A}N$ is defined as the completion of the quotient space of $M{\mathop{\widehat\otimes}} N$ over the closure of the linear span of all elements of the form
$$
x\cdot a\otimes y-x\otimes a\cdot y\qquad (x\in M,\,y\in N,\,a\in A).
$$
A left Fr\'echet module $F$ over a Fr\'echet algebra $A$ is said to be \emph{flat} if the functor $(-)\ptens{A} F$ maps every admissible complex of right Fr\'echet $A$-modules to an exact complex of vector spaces. It should be noted that we consider flat modules only in the category of Fr\'echet modules over a given Fr\'echet algebra (unlike projective modules); cf. Remark~\ref{onFR} below.

Recall that an \emph{$R$-${\mathop{\widehat\otimes}}$-algebra} is a pair $(A,\eta_A)$, where $A$ is a ${\mathop{\widehat\otimes}}$-algebra and $\eta_A \!: R \to A$ is a homomorphism of ${\mathop{\widehat\otimes}}$-algebras. Note that every $A$-${\mathop{\widehat\otimes}}$-module is automatically an $R$-${\mathop{\widehat\otimes}}$-module via the restriction of scalars functor along $\eta_A$. Denote by $(A,R){\mbox{-}\!\mathop{\mathsf{mod}}}$ the category of left $A$-${\mathop{\widehat\otimes}}$-modules endowed with the structure of an exact category with respect to complexes split by continuous morphisms of $R$-${\mathop{\widehat\otimes}}$-modules; cf. \cite[Appendix, Example 10.1 and 10.3]{Pir_qfree}. In particular, when $R = \mathbb{C}$, we obtain the standard definition of an admissible (or $\mathbb{C}$-split) sequence of $A$-${\mathop{\widehat\otimes}}$-modules mentioned above. Considering ${\mathop{\widehat\otimes}}$-bimodules over an $R$-${\mathop{\widehat\otimes}}$-algebra~$A$ (on the left) and an $S$-${\mathop{\widehat\otimes}}$-algebra $B$ (on the right), we denote the corresponding category by $(A,R){\mbox{-}\!\mathop{\mathsf{mod}}\!\mbox{-}}(B,S)$. In the case $R=S=\mathbb{C}$, we simply write $A{\mbox{-}\!\mathop{\mathsf{mod}}\!\mbox{-}} B$. For the categories of left and right modules, we use the notation $A{\mbox{-}\!\mathop{\mathsf{mod}}}$ and ${\mathop{\mathsf{mod}}\!\mbox{-}} B$, respectively. We also need the notion of relative projectivity. Namely, a module $P$ is called \emph{projective} in $(A,R){\mbox{-}\!\mathop{\mathsf{mod}}\!\mbox{-}}(B,S)$ if for any morphism with codomain $P$ having a right inverse morphism in $R{\mbox{-}\!\mathop{\mathsf{mod}}\!\mbox{-}} S$, there exists a right inverse morphism in $A{\mbox{-}\!\mathop{\mathsf{mod}}\!\mbox{-}} B$.

A \emph{projective (free, flat) resolution} of a $\mathbin{\widehat{\otimes}}$-module $M$ is an admissible complex $0\leftarrow M \leftarrow P_\bullet$ in which all modules $P_n$ ($n\ge 0$) are projective (free, flat). Every free resolution is projective and every projective resolution of a Fr\'echet module is flat.

Now we turn to the definitions of homological dimensions in the case of $\mathbin{\widehat{\otimes}}$-algebras. Let $A$ be a $\mathbin{\widehat{\otimes}}$-algebra and $M\in A{\mbox{-}\!\mathop{\mathsf{mod}}}$. Then the \emph{projective homological dimension} of $M$ is the length of its shortest projective resolution (denoted by $\mathop{\mathrm{dh}}_A M$). The \emph{global dimension} of~$A$ is defined as $$\mathop{\mathrm{dg}} A\!: = \sup\{\textstyle{\mathop{\mathrm{dh}}_A M}\!:\,M\in A{\mbox{-}\!\mathop{\mathsf{mod}}} \}.$$ The \emph{bidimension} of $A$ (denoted by $\mathop{\mathrm{db}} A$) is the projective homological dimension of~$A$ in $A{\mbox{-}\!\mathop{\mathsf{mod}}\!\mbox{-}} A$. See details in \cite{X1}.

The \emph{weak homological dimension} of a Fr\'echet module $M$ is the length of its shortest flat resolution in the category of Fr\'echet modules (denoted by $\mathop{\mathrm{w.dh}}_A M$). The \emph{weak global dimension} of~$A$ is defined as $$\mathop{\mathrm{w.dg}} A\!: = \sup\{\textstyle{\mathop{\mathrm{w.dh}}_A M}\!:\,M\in A{\mbox{-}\!\mathop{\mathsf{mod}}} \}.$$ The \emph{weak bidimension} of a Fr\'echet algebra $A$ (denoted by $\mathop{\mathrm{w.db}} A$) is the weak homological dimension of~$A$ in $A{\mbox{-}\!\mathop{\mathsf{mod}}\!\mbox{-}} A$. See details in \cite{Pi08}.

\subsection*{Properties of the bar-resolution}

By an \emph{augmented $\mathbin{\widehat{\otimes}}$-algebra} we mean a $\mathbin{\widehat{\otimes}}$-algebra~$A$ endowed with a homomorphism
$\varepsilon_A\!:A\to \mathbb{C}$. In this case, $\mathbb{C}$ is an $A$-$\mathbin{\widehat{\otimes}}$-bimodule with respect to the multiplications $a\cdot \lambda \!:= \varepsilon_A(a)\lambda$ and $\lambda \cdot a\!:=\lambda \varepsilon_A(a)$. Homomorphisms between augmented $\mathbin{\widehat{\otimes}}$-algebras are defined in the obvious way. Consider now the non-normalized bar-resolution
\begin{equation}\label{barA}
0\leftarrow \mathbb{C} \xleftarrow{\varepsilon_A} \beta_\bullet(A,\mathbb{C})
\end{equation}
of the module $\mathbb{C}$ in $A{\mbox{-}\!\mathop{\mathsf{mod}}}$; see, e.g., \cite[Chapter\,III, \S\,2.3]{X1}. Identifying $A\mathbin{\widehat{\otimes}} \mathbb{C}$ and~$A$, we may assume that $\beta_n(A,\mathbb{C})=A^{(n+1)\mathbin{\widehat{\otimes}}}$ and the differential $d_n\!: \beta_{n+1}(A,\mathbb{C})\to \beta_n(A,\mathbb{C})$ is given by the formula
\begin{multline*}
d_n(a_0\otimes\cdots \otimes a_{n+1})\!:=\\
\sum_{k=0}^n (-1)^k a_0\otimes\cdots \otimes
a_ka_{k+1}\otimes\cdots \otimes a_{n+1}+(-1)^{n+1}
a_0\otimes\cdots \otimes a_n \varepsilon_A(a_{n+1}).
\end{multline*}

Let now $H$ be a $\mathbin{\widehat{\otimes}}$-Hopf algebra and $A$ an $H$-$\mathbin{\widehat{\otimes}}$-module algebra. Under additional conditions, the complex \eqref{barA} (where $\varepsilon_A$ is the counit of $A$) can be regarded as a complex of $(A\mathop{\widehat{\#}} H)$-$\mathbin{\widehat{\otimes}}$-modules. Indeed, define an action of $H$ on $\beta_n(A,\mathbb{C})$ by the formula
$$
h\cdot(a_0\otimes\cdots \otimes a_{n})\!:=\sum  (h_{(0)}\cdot a_0)\otimes \cdots \otimes (h_{(n)}\cdot a_n).
$$
(Here we again use the generalized Sweedler notation.)
It is easy to verify that we obtain a structure of an $H$-$\mathbin{\widehat{\otimes}}$-module on $\beta_n(A,\mathbb{C})$.
To define a structure of an $(A\mathop{\widehat{\#}} H)$-$\mathbin{\widehat{\otimes}}$-module put $(a\mathop{\#}  h)\cdot x\!:=a\cdot(h\cdot x)$ for $a\in A$, $h\in H$ and $x$ in $\beta_n(A,\mathbb{C})$ or $\mathbb{C}$.

The statement of Part~(A) of the following proposition was proved in \cite[Lemma 3.2]{Pi4} for the case $n=0$, i.e., for $\beta_0(A,\mathbb{C})=A$.

\begin{pr}\label{beAHcom}
Let $H$ be a $\mathbin{\widehat{\otimes}}$-Hopf algebra and $A$ an $H$-$\mathbin{\widehat{\otimes}}$-module algebra.

\emph{(A)}~Then $\mathbb{C}$ and $\beta_n(A,\mathbb{C})$ are $(A\mathop{\widehat{\#}} H)$-$\mathbin{\widehat{\otimes}}$-modules \emph{(}$n\in\mathbb{Z}_+$\emph{)}.

Suppose, in addition, that $A$ is an $H$-$\mathbin{\widehat{\otimes}}$-module bialgebra. Then

\emph{(B)}~the counit $\varepsilon_A\!:A\to \mathbb{C}$ is a morphism of $(A\mathop{\widehat{\#}} H)$-$\mathbin{\widehat{\otimes}}$-modules, and the linear map
\begin{equation}\label{cAHAiso}
\mathbb{C}\ptens{A\mathop{\widehat{\#}} H}A\to \mathbb{C}\!: \lambda\otimes a\mapsto \lambda\varepsilon_A(a)
\end{equation}
is well defined and an isomorphism;

\emph{(C)}
the bar-resolution \eqref{barA} is a complex of $(A\mathop{\widehat{\#}} H)$-$\mathbin{\widehat{\otimes}}$-modules.
\end{pr}
Note that in the proof we only use the second equality from \eqref{coumo}. So the assumption that $A$ is an $H$-$\mathbin{\widehat{\otimes}}$-module bialgebra is slightly redundant.

\begin{proof}
(A)~The continuity of the multiplications is obvious. The associativity of the multiplication in $\mathbb{C}$ follows from~\eqref{exmusmpr}, the second equality in~\eqref{coumo}, and the formula $\sum \varepsilon_H(h_{(1)})\varepsilon_H( h_{(2)})=\varepsilon_H(h)$.

It suffices to verify the associativity of the multiplication in $\beta_n(A,\mathbb{C})$  for products of the form $(1\mathop{\#}  h)(a\mathop{\#} 1)$.
Let $a,a_0,\ldots a_n\in A$ and $h\in H$. Then
$$
(1\mathop{\#}  h)\cdot((a\mathop{\#} 1)\cdot(a_0\otimes\cdots \otimes a_{n}))=\sum  (h_{(0)}\cdot (aa_0))\otimes \cdots \otimes (h_{(n)}\cdot a_n).
$$
On the other hand,  we have $(1\mathop{\#}  h)(a\mathop{\#} 1)=\sum (h_{(-1)}\cdot a)\otimes h_{(0)}$ from \eqref{exmusmpr}. Therefore,
$$
((1\mathop{\#}  h)(a\mathop{\#} 1)\cdot(a_0\otimes\cdots \otimes a_{n}))=\sum  (h_{(-1)}\cdot a) (h_{(0)}\cdot a_0)\otimes \cdots \otimes (h_{(n)}\cdot a_n).
$$
By \eqref{Hmodalgcond}, the expressions on the right-hand side coincide. Hence, the multiplication is associative and thus
$\beta_n(A,\mathbb{C})$ is an $(A\mathop{\widehat{\#}} H)$-$\mathbin{\widehat{\otimes}}$-module.

(B)~Let $a,a'\in A$ and $h\in H$. Then
$$
\varepsilon_A((a\mathop{\#}  h)\cdot a')=\varepsilon_A(a(h\cdot a'))=\varepsilon_A(a)\,\varepsilon_A(h\cdot a')=\varepsilon_A(a)\,\varepsilon_H(h)\,\varepsilon_A(a').
$$
(The last equality holds by the second formula in~\eqref{coumo}.)

On the other hand,
$$(a\mathop{\#}  h)\cdot \varepsilon_A(a')=\varepsilon_A(a)\,\varepsilon_H(h)\,\varepsilon_A(a')$$
by the definition of the action of $A\mathop{\widehat{\#}} H$ on $\mathbb{C}$. This means that $\varepsilon_A$ is a $(A\mathop{\widehat{\#}} H)$-$\mathbin{\widehat{\otimes}}$-module morphism.
In particular, \eqref{cAHAiso} is well defined.

Note that $\mathbb{C}\ptens{A}A\to \mathbb{C}$ is an isomorphism. Since $A$ is a subalgebra of~$A\mathop{\widehat{\#}} H$, we have that $\mathbb{C}\ptens{A\mathop{\widehat{\#}} H}A$ is a quotient space of $\mathbb{C}\ptens{A}A$. Hence, \eqref{cAHAiso} is also an isomorphism.

(C)~By Part~(B), the counit $\varepsilon_A$ is a morphism of $(A\mathop{\widehat{\#}} H)$-$\mathbin{\widehat{\otimes}}$-modules. To verify that each of $d_n$ is a morphism of $(A\mathop{\widehat{\#}} H)$-$\mathbin{\widehat{\otimes}}$-modules, it suffices to check that it is a morphism of $H$-modules.

On the one hand,
\begin{multline*}
d_n(h\cdot(a_0\otimes\cdots \otimes a_{n+1}))=\\
=\sum_{k=0}^n (-1)^k  (h_{(0)}\cdot a_0)\otimes\cdots \otimes
 ((h_{(k)}\cdot a_k) (h_{(k+1)}\cdot a_{k+1}))\otimes\cdots \otimes (h_{(n+1)}\cdot a_{n+1})+\\
 +(-1)^{n+1}(h_{(0)}\cdot a_0)\otimes\cdots \otimes (h_{(n)}\cdot a_n)  \varepsilon_A(h_{(n+1)}\cdot a_{n+1}).
\end{multline*}
On the other hand,
\begin{multline*}
h\cdot(d_n(a_0\otimes\cdots \otimes a_{n+1}))=\\
=\sum_{k=0}^n (h_{(0)}\cdot a_0)\otimes\cdots \otimes
 (h_{(k)}\cdot (a_k a_{k+1}))\otimes\cdots \otimes (h_{(n+1)}\cdot a_{n+1})+\\
 +(-1)^{n+1}(h_{(0)}\cdot a_0)\otimes\cdots
\otimes (h_{(n)}\cdot (a_n \varepsilon_A(a_{n+1})).
\end{multline*}
Since $h_{(k)}\cdot (a_k a_{k+1})=\sum (h_{(k)}\cdot a_k) (h_{(k+1)}\cdot a_{k+1})$ by \eqref{Hmodalgcond}, it remains to verify the coincidence of the last factors in the last summands. Indeed, from~\eqref{coumo} we have
$$\varepsilon_A(h_{(n+1)}\cdot a_{n+1})=\varepsilon_H(h_{(n+1)})\varepsilon_A( a_{n+1}).$$
Thus,
$$\sum (h_{(n)}\cdot a_n)  \varepsilon_A(h_{(n+1)}\cdot a_{n+1})=\sum (h_{(n)}\cdot a_n) \varepsilon_H(h_{(n+1)})\varepsilon_A( a_{n+1}).$$
Further, it follows from the formula $h_{(n)}=\sum h_{(n)}\varepsilon_H(h_{(n+1)})$ that
$$h_{(n)}\cdot (a_n \varepsilon_A(a_{n+1})=\sum (h_{(n)}\cdot a_n)\varepsilon_H(h_{(n+1)})\varepsilon_A( a_{n+1}).$$
Thus, the last factors coincide, which completes the proof.
\end{proof}

\begin{co}\label{BShomdi0}
Let $H$ be a $\mathbin{\widehat{\otimes}}$-Hopf algebra and $A$ an $H$-$\mathbin{\widehat{\otimes}}$-module bialgebra. Then \eqref{barA} is a projective resolution of the module $\mathbb{C}$ in $A{\mbox{-}\!\mathop{\mathsf{mod}}}$ and also is  a complex of $(A\mathop{\widehat{\#}} H)$-$\mathbin{\widehat{\otimes}}$-modules.
\end{co}
\begin{proof}
The fact that \eqref{barA} is a projective resolution in $A{\mbox{-}\!\mathop{\mathsf{mod}}}$ is standard. On the other hand, it follows from Proposition~\ref{beAHcom}  that \eqref{barA} is a complex of $(A\mathop{\widehat{\#}} H)$-$\mathbin{\widehat{\otimes}}$-modules.
\end{proof}

\subsection*{Relative homology and dimensions}

Recall that an $R$-$\mathbin{\widehat{\otimes}}$-algebra $A$ is said to be \emph{relatively homologically trivial} if $A$ is a projective module in the relative category $(A,R){\mbox{-}\!\mathop{\mathsf{mod}}\!\mbox{-}}(A,R)$ \cite[Definition~4.8]{Ar_aahe}.

In what follows we need two simple lemmas. The first is derived directly from the definitions.
\begin{lm}\label{SBproj}
If a module is projective in $(A,R){\mbox{-}\!\mathop{\mathsf{mod}}}$ and $R{\mbox{-}\!\mathop{\mathsf{mod}}}$, then it is projective in $A{\mbox{-}\!\mathop{\mathsf{mod}}}$.
\end{lm}

\begin{lm}\label{hotrallpr}
If an $R$-$\mathbin{\widehat{\otimes}}$-algebra $A$ is relatively homologically trivial, then all left $A$-$\mathbin{\widehat{\otimes}}$-modules are projective in $(A,R){\mbox{-}\!\mathop{\mathsf{mod}}}$.
\end{lm}
\begin{proof}
The argument is standard. Indeed, let $N$ be a left $A$-$\mathbin{\widehat{\otimes}}$-module. It suffices to show that the multiplication morphism $A\ptens{R}N\to N$ has a right inverse morphism in $A{\mbox{-}\!\mathop{\mathsf{mod}}}$. Note that the formula $a\mapsto 1\otimes a$ gives a right inverse in $R{\mbox{-}\!\mathop{\mathsf{mod}}\!\mbox{-}} R$ to the multiplication morphism $\pi\!:A\ptens{R}A\to A$. Since $A$ is a projective module in $(A,R){\mbox{-}\!\mathop{\mathsf{mod}}\!\mbox{-}}(A,R)$, we conclude that $\pi$ has a right inverse in $A{\mbox{-}\!\mathop{\mathsf{mod}}\!\mbox{-}} A$. Applying the functor $(-)\ptens{A} N$ to it and using the isomorphism $A\ptens{A} N\cong N$, we obtain the required right inverse in $A{\mbox{-}\!\mathop{\mathsf{mod}}}$.
\end{proof}

\begin{pr}\label{BShomdi}
Let $N$ be a left $\mathbin{\widehat{\otimes}}$-module over an $R$-$\mathbin{\widehat{\otimes}}$-algebra $A$. Suppose that $A$ is relatively homologically trivial. If there is a projective resolution of $N$ in $R{\mbox{-}\!\mathop{\mathsf{mod}}}$ that is also a complex in $A{\mbox{-}\!\mathop{\mathsf{mod}}}$, and $\mathop{\mathrm{dh}}_R N\le k$ for some $k\in\mathbb{Z}_+$, then $\mathop{\mathrm{dh}}_A N\le k$.
\end{pr}

\begin{proof}
Let $0\leftarrow N \leftarrow P_\bullet$ be a projective resolution of $N$ in $R{\mbox{-}\!\mathop{\mathsf{mod}}}$ that is also a complex in $A{\mbox{-}\!\mathop{\mathsf{mod}}}$.
Consider the complex
\begin{equation}\label{cowker}
0\leftarrow N \leftarrow P_0\xleftarrow{d_0}\cdots\xleftarrow{d_{k-2}} P_{k-1} \leftarrow \Ker d_{k-2}\leftarrow 0.
\end{equation}
Since it is admissible, by a standard result of homological algebra, the inequality $\mathop{\mathrm{dh}}_R N\le k$ implies that not only $P_0,\ldots,P_{k-1}$, but also $\Ker d_{k-2}$ are projective in $R{\mbox{-}\!\mathop{\mathsf{mod}}}$. On the other hand,  all modules are projective in $(A,R){\mbox{-}\!\mathop{\mathsf{mod}}}$ by Lemma~\ref{hotrallpr}. Then we obtain from Lemma~\ref{SBproj} that $P_0,\ldots,P_{k-1}$ and $\Ker d_{k-2}$ are also projective in $A{\mbox{-}\!\mathop{\mathsf{mod}}}$. Since~\eqref{cowker} is a complex in $A{\mbox{-}\!\mathop{\mathsf{mod}}}$, it is a projective resolution of length $k$ in $A{\mbox{-}\!\mathop{\mathsf{mod}}}$. Thus $\mathop{\mathrm{dh}}_A N\le k$.
\end{proof}

\subsection*{Homological epimorphisms}

A $\mathbin{\widehat{\otimes}}$-algebra homomorphism  $\varphi\!:A\to B$  is called a \emph{homological epimorphism} if for some (or, equivalently, for every) projective resolution $0 \leftarrow B \leftarrow P_\bullet$ in $A{\mbox{-}\!\mathop{\mathsf{mod}}}$, the sequence
$$
0 \longleftarrow B \longleftarrow B\ptens{A} P_\bullet
$$
is admissible; see, e.g., \cite[Definition 6.3]{Pir_stbflat}. (This notion may appear under different names; see \cite[Remark 3.16]{AP}.)

Recall that for $M\in {\mathop{\mathsf{mod}}\!\mbox{-}} A$, $N\in A{\mbox{-}\!\mathop{\mathsf{mod}}}$ and $n\in\mathbb{Z}_+$, the derived (in the classical sense) functor ${\mathop{\mathrm{Tor}}\nolimits}_n^{A}(M,N)$  can be defined as the homology of the complex $M\ptens{A} P_\bullet$ (in the category of vector spaces), where $P_\bullet$ is an arbitrary projective resolution of $N$; see \cite[Chapter III, \S\,4]{X1} or \cite{He81}.
Following \cite[Appendix, p.\,168, Definition~1]{He81}, we say that a  $\mathbin{\widehat{\otimes}}$-algebra homomorphism  $\varphi\!: A \to B$ is \emph{torsion-preserving} if all linear maps
$$
{\mathop{\mathrm{Tor}}\nolimits}_n^{A}(M,N)\to{\mathop{\mathrm{Tor}}\nolimits}_n^{B}(M,N) \qquad (M\in {\mathop{\mathsf{mod}}\!\mbox{-}} B,\, N\in B{\mbox{-}\!\mathop{\mathsf{mod}}},\,n\in\mathbb{Z}_+),
$$
induced by it, are bijective.

\begin{pr}\label{Hptor}
\emph{(Taylor)}
Every homological epimorphism of $\mathbin{\widehat{\otimes}}$-algebras is torsion-preserving.
\end{pr}
\begin{proof}
According to \cite[Proposition 1.6]{T2}, every homological epimorphism is a Taylor localization with respect to the class of all complete locally convex spaces. By \cite[Proposition 1.4]{T2}, every Taylor localization is torsion-preserving.
\end{proof}

The following lemma is needed for the proof of Theorem~\ref{toriAG}.
\begin{lm}\label{quprTor}
Let $\varphi\!:A\to B$ be a left (or right) projective epimorphism of $\mathbin{\widehat{\otimes}}$-algebras, i.e., $\varphi$ is an epimorphism and $B$ is projective as a left (or right) $A$-$\mathbin{\widehat{\otimes}}$-module. Then $\varphi$ is a homological epimorphism.
\end{lm}
\begin{proof}
The following simple fact is well known (see, e.g., \cite[Proposition 6.1]{Pir_qfree}):
if $\varphi$ is an epimorphism, then $B\ptens{A}B\to B$ is a topological isomorphism of modules.

Let $0 \leftarrow B\leftarrow Q_\bullet$ be a free resolution in $B{\mbox{-}\!\mathop{\mathsf{mod}}}$. If $B$ is a projective left $A$-$\mathbin{\widehat{\otimes}}$-module, then this resolution is an admissible complex in $A{\mbox{-}\!\mathop{\mathsf{mod}}}$, each term of which is projective. From the above property of epimorphisms, it follows that the functor $B\ptens{A}(-)$ maps this resolution to itself. In particular, it preserves its admissibility. This means that $\varphi$ is a homological epimorphism.
\end{proof}

\subsection*{Theorem on the completion of ${\mathscr A}(G)$}

For an arbitrary complex Lie group $G$, denote by $\Lambda$ its linearizer, i.e., the intersection of the kernels of all finite-dimensional holomorphic representations. In many questions, it is reasonable to replace $G$ by its linearization $G/\Lambda$, which is automatically linear.
In the case when $G$ is connected and linear, one can define its nilpotent and exponential radicals. The former is defined as the intersection of the kernels of all \emph{irreducible} finite-dimensional holomorphic representations, and for the definition of the latter we refer to \cite[\S\,3]{ArAnF}. We also omit the definition of a submultiplicative weight with exponential distortion on a given
subgroup \cite[Definition 4.1]{Ar_smash}. It is essential for us that there exists a maximal (up to the natural equivalence relation) weight among such weights for a fixed (normal integral) subgroup intermediate between the two radicals \cite[Theorem 4.3]{Ar_smash}; cf. \cite[Definition 2.4]{ArLfd}. Denote this weight by $\omega_{max}$. Since $\omega_{max}$ is submultiplicative, one can associate with it a Fr\'echet--Arens--Michael algebra, which is a completion of the algebra ${\mathscr A}(G/\Lambda)$ and which we denote by ${\mathscr A}_{\omega_{max}^\infty}(G/\Lambda)$; for details see~\cite[\S\,3]{Ar_smash}.

In estimating the homological dimensions of ${\mathscr A}_{\omega_{max}^\infty}(G/\Lambda)$, we  rely on the following result.
\begin{thm}\label{C4mainexdi}  \cite[Theorem 3.2]{Ar_aahe}
Let $G$ be a connected complex Lie group, $\Lambda$ its linearizer, and $E$ and $N$ the exponential and nilpotent radicals of~$G/\Lambda$. Suppose that $N'$ is a normal integral subgroup of $G/\Lambda$ such that $E\subset N'\subset N$ and $\omega_{max}$ is a submultiplicative weight on~$G/\Lambda$ that is maximal among weights with exponential distortion on $N'$.
Then ${\mathscr A}(G)\to {\mathscr A}_{\omega_{max}^\infty}(G/\Lambda)$ is a homological epimorphism.
\end{thm}

\section{Computation of homological dimensions}
\label{s:homdim}

Now we formulate the central result of the article. Recall that a connected linear complex Lie group can be represented as $B\rtimes L$, where $B$ is simply connected and solvable, and $L$ is linearly complex reductive \cite[p.\,601, Theorem 16.3.7]{HiNe}. In particular, we may apply this result for $G/\Lambda$ if $G$ is an arbitrary connected Lie group and $\Lambda$ is its linearizer. Now fix a decomposition $G/\Lambda=B\rtimes L$.

Let $E$ and $N$ be the exponential and nilpotent radicals of~$G/\Lambda$, and let $N'$ be a normal integral subgroup of $G/\Lambda$ such that $E\subset N'\subset N$. Denote, as above, the submultiplicative weight on~$G/\Lambda$ that is maximal among weights with exponential distortion on $N'$ by $\omega_{max}$, and additionally denote the algebra ${\mathscr A}_{\omega_{max}^\infty}(G/\Lambda)$ by~$D$.

\begin{thm}\label{homdimth}
Let $G$ be a connected complex Lie group. Then
$$
\mathop{\mathrm{dg}} D=\mathop{\mathrm{dg}} {\mathscr A}(G)=\mathop{\mathrm{db}} D=\mathop{\mathrm{db}} {\mathscr A}(G)={\textstyle
\mathop{\mathrm{dh}}_D}\mathbb{C}={\textstyle\mathop{\mathrm{dh}}_{{\mathscr A}(G)}}\mathbb{C}=\dim B.
$$
Similarly, for weak dimensions we have
$$
\mathop{\mathrm{w.dg}} D=\mathop{\mathrm{w.db}} D={\textstyle\mathop{\mathrm{w.dh}}_D}\mathbb{C}=\dim B.
$$
\end{thm}

Note that the most important step in the argument linking the homological properties of ${\mathscr A}(G)$ and its completion $D$ --- the existence of the isomorphism ${\mathop{\mathrm{Tor}}\nolimits}_n^{D}(\mathbb{C},\mathbb{C})\cong{\mathop{\mathrm{Tor}}\nolimits}_n^{{\mathscr A}(G)}(\mathbb{C},\mathbb{C})$ in Theorem~\ref{toriAG} --- uses the results on homological epimorphisms obtained in \cite{Ar_aahe}, in particular Theorem~\ref{C4mainexdi} formulated above.

\begin{rem}\label{onFR}
Theorem~\ref{homdimth} does not include the computation of weak dimensions for ${\mathscr A}(G)$. The reason is that for $\mathbin{\widehat{\otimes}}$-algebras that are not Fr\'echet spaces, the very notion of a flat module (and hence of weak dimensions) is ambiguous; see \cite[Definition~3.1, Remark 3.2]{PP22}. However, the author has no doubt that with a properly chosen definition, analogous formulas should also hold for the weak dimensions of ${\mathscr A}(G)$.
\end{rem}

Since the subgroup $N'$ in Theorem~\ref{homdimth} can be chosen such that ${\mathscr A}_{\omega_{max}^\infty}(G/\Lambda)\cong \widehat{\mathscr A}(G)$, or such that ${\mathscr A}_{\omega_{max}^\infty}(G/\Lambda)\cong \widehat{\mathscr A}(G)^{\mathrm{PI}}$ \cite[Theorem 5.11]{Ar_smash}, we obtain two corollaries for these special cases.

\begin{co}\label{hdAM}
Let $G$ be a connected complex Lie group. Then
 $$
\mathop{\mathrm{dg}} \widehat{\mathscr A}(G)=\mathop{\mathrm{db}} \widehat{\mathscr A}(G)={\textstyle
\mathop{\mathrm{dh}}_{\widehat{\mathscr A}(G)}}\mathbb{C}=\dim B;
$$
$$
\mathop{\mathrm{w.dg}} \widehat{\mathscr A}(G)=\mathop{\mathrm{w.db}} \widehat{\mathscr A}(G)={\textstyle\mathop{\mathrm{w.dh}}_{\widehat{\mathscr A}(G)}}\mathbb{C}=\dim B.
$$
\end{co}

\begin{co}\label{hdPIE}
Let $G$ be a connected complex Lie group. Then
 $$
\mathop{\mathrm{dg}} \widehat{\mathscr A}(G)^{\mathrm{PI}}=\mathop{\mathrm{db}} \widehat{\mathscr A}(G)^{\mathrm{PI}}={\textstyle
\mathop{\mathrm{dh}}_{\widehat{\mathscr A}(G)^{\mathrm{PI}}}}\mathbb{C}=\dim B;
$$
 $$
\mathop{\mathrm{w.dg}} \widehat{\mathscr A}(G)^{\mathrm{PI}}=\mathop{\mathrm{w.db}} \widehat{\mathscr A}(G)^{\mathrm{PI}}={\textstyle
\mathop{\mathrm{w.dh}}_{\widehat{\mathscr A}(G)^{\mathrm{PI}}}}\mathbb{C}=\dim B.
$$
\end{co}

For the proof of Theorem~\ref{homdimth}, we need some auxiliary facts. First, we use a statement about the bar-resolution~\eqref{barA}.

\begin{pr}\label{beAB}
The complex
\begin{equation}\label{ABCC}
0\leftarrow \mathbb{C} \xleftarrow{\varepsilon_{{\mathscr A}(B)}} \beta_\bullet({\mathscr A}(B),\mathbb{C})
\end{equation}
is a projective resolution in ${\mathscr A}(G/\Lambda){\mbox{-}\!\mathop{\mathsf{mod}}}$.
\end{pr}
\begin{proof}
Consider the semidirect product $B\rtimes L$. It follows from \cite[Theorem 2.1]{Ar_smash} that ${\mathscr A}(B)$ is an ${\mathscr A}(L)$-$\mathbin{\widehat{\otimes}}$-module bialgebra. Then by Corollary~\ref{BShomdi0}, the sequence~\eqref{ABCC} is a complex of $({\mathscr A}(B)\mathop{\widehat{\#}}{\mathscr A}(L))$-$\mathbin{\widehat{\otimes}}$-modules.

Since $L$ is linearly complex reductive, ${\mathscr A}(L)$ is homologically trivial \cite[Theorem 4.1]{Ar_aahe}. Thus, all left ${\mathscr A}(L)$-$\mathbin{\widehat{\otimes}}$-modules are projective, in particular~$\mathbb{C}$. Then ${\mathscr A}(B)$ is a projective left $({\mathscr A}(B)\mathop{\widehat{\#}}{\mathscr A}(L))$-$\mathbin{\widehat{\otimes}}$-module by \cite[Lemma~3.3]{Pi4}. As a consequence, the module $\beta_n({\mathscr A}(B),\mathbb{C})$ is projective in $({\mathscr A}(B)\mathop{\widehat{\#}}{\mathscr A}(L)){\mbox{-}\!\mathop{\mathsf{mod}}}$ for each $n$. The aforementioned Theorem~2.1 in \cite{Ar_smash} implies that ${\mathscr A}(B)\mathop{\widehat{\#}} {\mathscr A}(L)\cong {\mathscr A}(B\rtimes L)$. Recalling that
${\mathscr A}(B\rtimes L)\cong {\mathscr A}(G/\Lambda)$ completes the proof.
\end{proof}

Second, let $G$ be a complex Lie group with the Lie algebra~$\mathfrak{g}$, the universal enveloping algebra $U(\mathfrak{g})$, and the standard embedding $\tau\!:U(\mathfrak{g}) \to {\mathscr A}(G)$. Recall that $U(\mathfrak{g})$ is a right $\mathfrak{g}$-module with respect to the right regular representation given by the formula $a\cdot X\!:= aX$, where $a\in U(\mathfrak{g})$ and $X\in\mathfrak{g}$. Similarly, ${\mathscr A}(G)$ is a right $\mathfrak{g}$-module with respect to the multiplication $\mu\cdot X\!:= \mu\tau(X)$.

Let $C_\bullet(\mathfrak{g},M)$ be the standard chain complex for a right $\mathfrak{g}$-module~$M$ (called in this context the Koszul complex or the Chevalley--Eilenberg chain complex), i.e.,
\begin{equation}\label{CCEco}
0\leftarrow M\leftarrow M\otimes \textstyle{\bigwedge^1}\mathfrak{g}
\leftarrow \cdots \leftarrow  M\otimes\textstyle{\bigwedge^n}\mathfrak{g} \leftarrow
\cdots   \,,
\end{equation}
where the differential $\delta_n\!:M\otimes \textstyle{\bigwedge^{n+1}}\mathfrak{g}\to
M\otimes\textstyle{\bigwedge^n}\mathfrak{g}$ is given by the formula
\begin{multline*}
\delta_n(m\otimes X_1\wedge\cdots\wedge X_{n+1})= \sum_{i=1}^{n+1}
(-1)^{i-1} m\cdot X_i\otimes X_1\wedge\cdots\wedge\hat X_i\wedge
\cdots\wedge X_{n+1}\\
+\sum_{1\le i<j\le {n+1}} (-1)^{i+j} m\otimes [X_i,X_j]\wedge
X_1\wedge\cdots\wedge\hat X_i \wedge\cdots\wedge\hat
X_j\wedge\cdots\wedge X_{n+1}.
\end{multline*}
(Here the notation $\hat X_i$ means that $X_i$ is omitted.) The $n$th homology of this complex is denoted by $\mathrm{H}_n^{\mathrm{Lie}}(\mathfrak{b},M)$ and is nothing other than the Lie algebra homology of $\mathfrak{g}$ with coefficients in~$M$.
A classical theorem (see, e.g., \cite[p.\,339, Chapter XIII, Theorem 7.1]{CE60}) states that the complex
\begin{equation}\label{CUgres}
0 \leftarrow\mathbb{C}\xleftarrow{\varepsilon}  C_\bullet(\mathfrak{g},U(\mathfrak{g})),
\end{equation}
where $\varepsilon: U(\mathfrak{g})\to\mathbb{C}$ is the counit, is exact.

The following result, which is also of independent interest, is needed. We use the notation $\Lambda$, $L$, $B$, and $D$ as before. Note that the chains of homomorphisms
$$
 {\mathscr A}(G)\to {\mathscr A}(G/\Lambda)\to D \quad\text{and}\quad U(\mathfrak{b})\to {\mathscr A}(B)\to {\mathscr A}(G/\Lambda)
$$
transform every $D$-$\mathbin{\widehat{\otimes}}$-module into a $\mathbin{\widehat{\otimes}}$-module over the remaining algebras.

\begin{thm}\label{toriAG}
Let $G$ be a connected complex Lie group and $M$ a one-dimensional right $D$-$\mathbin{\widehat{\otimes}}$-module. Then for each $n\in\mathbb{Z}_+$
$$
{\mathop{\mathrm{Tor}}\nolimits}_n^{D}(M,\mathbb{C})\cong{\mathop{\mathrm{Tor}}\nolimits}_n^{{\mathscr A}(G)}(M,\mathbb{C})\cong{\mathop{\mathrm{Tor}}\nolimits}_n^{{\mathscr A}(G/\Lambda)}(M,\mathbb{C})\cong
{\mathop{\mathrm{Tor}}\nolimits}_n^{U(\mathfrak{b})}(M,\mathbb{C})\cong \textstyle{\mathrm{H}_n^{\mathrm{Lie}}(\mathfrak{b},M)},
$$
where $\mathfrak{b}$ is the Lie algebra associated with~$B$. (Here vector space isomorphisms are implied.)
\end{thm}

Note that the algebra $U(\mathfrak{b})$ has countable linear dimension and so it can be treated as a $\mathbin{\widehat{\otimes}}$-algebra with respect to the strongest locally convex topology. Therefore, it is irrelevant here in which sense the functor ${\mathop{\mathrm{Tor}}\nolimits}_n^{U(\mathfrak{b})}(-,-)$ is understood, purely algebraically or taking topology into account.

Before proving the theorem, we construct a resolution. Take the complex~\eqref{CUgres} with~$\mathfrak{b}$ instead of~$\mathfrak{g}$, i.e.,
\begin{equation}\label{CUbres}
0 \leftarrow\mathbb{C}\xleftarrow{\varepsilon}  C_\bullet(\mathfrak{b},U(\mathfrak{b})).
\end{equation}
As noted above, it is exact and $U(\mathfrak{b})$ can be endowed with the strongest locally convex topology. Thus, \eqref{CUbres} is a projective resolution in $U(\mathfrak{b}){\mbox{-}\!\mathop{\mathsf{mod}}}$.
Since the image of $U(\mathfrak{b})$ is dense in ${\mathscr A}(B)$, we have
$$
{\mathscr A}(B)\ptens{U(\mathfrak{b})}\mathbb{C}\cong {\mathscr A}(B)\ptens{{\mathscr A}(B)}\mathbb{C}\cong \mathbb{C}
$$
in ${\mathscr A}(B){\mbox{-}\!\mathop{\mathsf{mod}}}$. It is easy to check that the standard embedding $U(\mathfrak{b})\to {\mathscr A}(B)$ is a homomorphism of augmented $\mathbin{\widehat{\otimes}}$-algebras.
From this it is easy to see that applying the functor ${\mathscr A}(B)\ptens{U(\mathfrak{b})}(-)$ to~\eqref{CUbres}, we obtain a complex in ${\mathscr A}(B){\mbox{-}\!\mathop{\mathsf{mod}}}$, isomorphic to
\begin{equation}\label{CbbAB}
0 \leftarrow\mathbb{C}\leftarrow  C_\bullet(\mathfrak{b},{\mathscr A}(B)),
\end{equation}
where the morphism ${\mathscr A}(B)\to\mathbb{C}$ coincides with the counit.

\begin{lm}\label{GLaBres}
The complex~\eqref{CbbAB} is a projective resolution in ${\mathscr A}(B){\mbox{-}\!\mathop{\mathsf{mod}}}$ and, consequently, $\mathop{\mathrm{dh}}_{{\mathscr A}(B)}\mathbb{C}\le k$.
\end{lm}
\begin{proof}
Since $B$ is simply connected and solvable, $\tau\!:U(\mathfrak{b})\to {\mathscr A}(B)$ is a homological epimorphism \cite[Theorem~8.3]{Pir_stbflat}. Consequently, by \cite[Proposition 3.2]{Pir_stbflat}, the
application the functor ${\mathscr A}(B)\ptens{U(\mathfrak{b})}(-)$ to a projective resolution of $\mathbb{C}$ in $U(\mathfrak{b}){\mbox{-}\!\mathop{\mathsf{mod}}}$ yields an admissible complex. Thus, \eqref{CbbAB} is admissible and is therefore a projective resolution in ${\mathscr A}(B){\mbox{-}\!\mathop{\mathsf{mod}}}$. Since~$\mathfrak{b}$ has linear dimension $k$, we conclude that ${\bigwedge^n}\mathfrak{g}=0$ for $n>k$. Hence the resolution has length~$k$.
\end{proof}

\begin{proof}[Proof of Theorem~\ref{toriAG}]
(1)~By Theorem~\ref{C4mainexdi}, the homomorphism ${{\mathscr A}}(G)\to D$ is a homological epimorphism. Also by Proposition~\ref{Hptor}, each homological epimorphism is torsion-preserving. In particular, ${\mathop{\mathrm{Tor}}\nolimits}_n^{D}(M,\mathbb{C})\cong{\mathop{\mathrm{Tor}}\nolimits}_n^{{\mathscr A}(G)}(M,\mathbb{C})$ for every~$n$.

(2)~We now claim that ${\mathop{\mathrm{Tor}}\nolimits}_n^{{\mathscr A}(G)}(M,\mathbb{C})\cong{\mathop{\mathrm{Tor}}\nolimits}_n^{{\mathscr A}(G/\Lambda)}(M,\mathbb{C})$ for every~$n$. Indeed, the $\mathbin{\widehat{\otimes}}$-algebra ${\mathscr A}(\Lambda)$ (where, as above, $\Lambda$ is the linearizer of~$G$) is homologically trivial \cite[Corollary~5.4]{Ar_aahe}. Then \cite[Proposition~5.7]{Ar_aahe} implies that ${\mathscr A}(G/\Lambda)$ is a projective left
${\mathscr A}(G)$-$\mathbin{\widehat{\otimes}}$-module. Moreover, the homomorphism ${\mathscr A}(G)\to{\mathscr A}(G/\Lambda)$ is surjective
\cite[Proposition~3.9]{AHHFG} and so it is torsion-preserving by Lemma~\ref{quprTor}. In particular, ${\mathop{\mathrm{Tor}}\nolimits}_n^{{\mathscr A}(G)}(M,\mathbb{C})\cong{\mathop{\mathrm{Tor}}\nolimits}_n^{{\mathscr A}(G/\Lambda)}(M,\mathbb{C})$ for every $n$.

(3)~To complete the proof we show that ${\mathop{\mathrm{Tor}}\nolimits}_n^{{\mathscr A}(G/\Lambda)}(M,\mathbb{C})\cong{\mathop{\mathrm{Tor}}\nolimits}_n^{U(\mathfrak{b})}(M,\mathbb{C})\cong \mathrm{H}^{\mathrm{Lie}}_n(\mathfrak{b},M)$ for every~$n$. For this, we use two resolutions. First, we claim that the map ${\mathop{\mathrm{Tor}}\nolimits}_n^{{\mathscr A}(G/\Lambda)}(M,\mathbb{C})\to {\mathop{\mathrm{Tor}}\nolimits}_n^{{\mathscr A}(B)}(M,\mathbb{C})$, induced by the embedding $B\to G/\Lambda$, is a bijection. Indeed, by Proposition~\ref{beAB}, the complex $0\leftarrow \mathbb{C} \leftarrow  \beta_\bullet({\mathscr A}(B),\mathbb{C})$ is a projective resolution in ${\mathscr A}(G/\Lambda){\mbox{-}\!\mathop{\mathsf{mod}}}$. Moreover, it is obviously a projective resolution in ${\mathscr A}(B){\mbox{-}\!\mathop{\mathsf{mod}}}$. Thus, we can use it to compute both ${\mathop{\mathrm{Tor}}\nolimits}_n^{{\mathscr A}(G/\Lambda)}(M,\mathbb{C})$ and ${\mathop{\mathrm{Tor}}\nolimits}_n^{{\mathscr A}(B)}(M,\mathbb{C})$.

Obviously, $M\ptens{{\mathscr A}(B)}{\mathscr A}(B)\cong M$. Furthermore, we claim that $M\ptens{{\mathscr A}(G/\Lambda)}{\mathscr A}(B)\cong M$. Indeed, since  $M$ is one-dimensional, it suffices to show that $M\ptens{{\mathscr A}(G/\Lambda)}{\mathscr A}(B)\ne 0$. Consider the dual left ${\mathscr A}(G/\Lambda)$-module $M^*$. Since it is a quotient module of ${\mathscr A}(B)$, there exists a continuous linear map $$M\ptens{{\mathscr A}(G/\Lambda)}{\mathscr A}(B)\to M\ptens{{\mathscr A}(G/\Lambda)}M^*.$$ The linear functional $M{\mathop{\widehat\otimes}} M^*\to \mathbb{C}$ associated with the natural bilinear functional $M\times M^*\to \mathbb{C}$ is not degenerate and factors through $M\ptens{{\mathscr A}(G/\Lambda)}M^*$. Consequently, $M\ptens{{\mathscr A}(G/\Lambda)}M^*$ and hence $M\ptens{{\mathscr A}(G/\Lambda)}{\mathscr A}(B)$ are both non-trivial.

Further, it is easy to see that application of the functors $M\ptens{{\mathscr A}(B)}(-)$ and $M\ptens{{\mathscr A}(G/\Lambda)}(-)$ to the complex $\beta_\bullet({\mathscr A}(B),\mathbb{C})$ gives the same result. Thus, ${\mathop{\mathrm{Tor}}\nolimits}_n^{{\mathscr A}(G/\Lambda)}(M,\mathbb{C})\cong {\mathop{\mathrm{Tor}}\nolimits}_n^{{\mathscr A}(B)}(M,\mathbb{C})$.

Next, we claim that ${\mathop{\mathrm{Tor}}\nolimits}_n^{{\mathscr A}(B)}(M,\mathbb{C})\cong \mathrm{H}^{\mathrm{Lie}}_n(\mathfrak{b},M)$. Indeed, by Lemma~\ref{GLaBres}, we can use the resolution~\eqref{CbbAB} to compute ${\mathop{\mathrm{Tor}}\nolimits}_n^{{\mathscr A}(B)}(M,\mathbb{C})$. Applying the functor $M\ptens{{\mathscr A}(B)}(-)$ to $C_\bullet(\mathfrak{b},{\mathscr A}(B))$, we obtain the complex $C_\bullet(\mathfrak{b},M)$, whose homology coincides with $\mathrm{H}^{\mathrm{Lie}}_n(\mathfrak{b},\mathbb{C})$ by definition.

Finally, note that the existence of the isomorphism ${\mathop{\mathrm{Tor}}\nolimits}_n^{U(\mathfrak{b})}(M,\mathbb{C})\cong \mathrm{H}^{\mathrm{Lie}}_n(\mathfrak{b},M)$ for each $n\in\mathbb{Z}_+$ is a standard fact from the theory of Lie algebra homology.
\end{proof}

Since $B$ is a normal subgroup of $G/\Lambda$, we can take the restriction of the adjoint representation of $G/\Lambda$ to $\mathfrak{b}$, which we denote by ${\mathop{\mathrm{Ad}}\nolimits}_\mathfrak{b}$.
Consider also the one-dimensional representation of $G/\Lambda$ given by the formula
\begin{equation}\label{ADbdet}
\chi\!: h\mapsto \det({\mathop{\mathrm{Ad}}\nolimits}_\mathfrak{b} h)
\end{equation}
and denote the corresponding left module by $\mathbb{C}_\chi$. Note that $\chi$ extends to ${\mathscr A}(G/\Lambda)$, being a holomorphic representation of $G/\Lambda$.

\begin{lm}\label{prochi}
\emph{(A)}~The representation $\chi$ extends to $D$.

\emph{(B)}~Being restricted along $\mathfrak{b}\to {\mathscr A}(G/\Lambda)$, the representation $\chi$ takes the form $X\mapsto \mathop{\mathrm{Tr}}\nolimits({\mathop{\mathrm{ad}}\nolimits} X)$.
\end{lm}
\begin{proof}
(A)~Since $N$ is the intersection of the kernels of all irreducible finite-dimensional holomorphic representations, we have $\chi(N)=1$. Then since $N'\subset N$, we conclude that $g\mapsto |\chi(g)|$ satisfies the inequality in \cite[Definition 4.1]{Ar_smash} and hence has exponential distortion on $N'$. Recall that by definition $D={\mathscr A}_{\omega_{max}^\infty}(G/\Lambda)$, where $\omega_{max}$ is a submultiplicative weight on~$G/\Lambda$ that is maximal among weights with exponential distortion on~$N'$. Therefore, by \cite[Theorem 5.1]{Ar_smash}, the algebra $D$ has a universal property for holomorphic homomorphisms whose norms have exponential distortion. It follows from the indicated universal property that~$\chi$ extends to~$D$.

(B)~It suffices to prove the statement for the restriction of $\chi$ to~$B$. Recall that for a holomorphic matrix-valued function $z\mapsto A(z)$, the equality
$$
\frac{d(\det A(z))}{dz}= \det A(z)\,\mathop{\mathrm{Tr}}\nolimits\left(\frac{dA(z)}{dz} A(z)^{-1}\right)
$$
holds; see, e.g., \cite[Chapter~8, p.\,443, Problem 8.6]{Pr08}. For $X\in \mathfrak{b}$, putting $A(z)={\mathop{\mathrm{Ad}}\nolimits} \exp(zX)$ and taking into account that $A(0)=1$, we obtain the equality
$$
\frac{d}{dz}\Bigl |_{z=0}\chi(\exp(zX))=\mathop{\mathrm{Tr}}\nolimits\frac{d({\mathop{\mathrm{Ad}}\nolimits} \exp(zX))}{dz}\Bigl |_{z=0}.
$$
As is well known, the right-hand side is equal to ${\mathop{\mathrm{ad}}\nolimits} X$. On the other hand, the discussed representation of the algebra $U(\mathfrak{b})$ is generated by the map sending $X$ to the left-hand side of the equality. This proves Part~(B).
\end{proof}

Now everything is ready for the proof of the main result. For the upper bound, we use homological epimorphisms, and for the lower bound, the result on the nontriviality of Lie algebra cohomology from \cite{Ha70}.

\begin{proof}[Proof of Theorem~\ref{homdimth}]
According to \cite[Corollary 2.7]{Pir_stbflat},
$$
\textstyle{\mathop{\mathrm{dh}}_H}\mathbb{C}=\mathop{\mathrm{dg}} H=\mathop{\mathrm{db}} H
$$
for every $\mathbin{\widehat{\otimes}}$-Hopf algebra $H$ with invertible antipode. It is easy to show that for every cocommutative $\mathbin{\widehat{\otimes}}$-Hopf algebra, the square of the antipode coincides with the identity map; see \cite[Proposition~5.9]{Ar_aahe}. Thus, we can use this formula in the cases where $H$ is ${\mathscr A}(G)$ or $D$. Furthermore,
$$
\textstyle{\mathop{\mathrm{w.dh}}_D}\mathbb{C}\le\mathop{\mathrm{w.dg}} D\le\mathop{\mathrm{w.db}} D,
$$
where the first inequality is an obvious consequence of the definition of weak global dimension and the second was proved in \cite[Proposition 4.7]{Pi08} for Fr\'echet algebras. Furthermore, from the fact that every projective module is flat, we obtain the inequalities
$$
\textstyle{\mathop{\mathrm{w.dh}}_D}\mathbb{C}\le \textstyle{\mathop{\mathrm{dh}}_D}\mathbb{C},\quad \mathop{\mathrm{w.dg}} D\le\mathop{\mathrm{dg}} D,\quad \mathop{\mathrm{w.db}} D\le \mathop{\mathrm{db}} D.
$$

Thus, it suffices to establish that
$$
k\le\textstyle{\mathop{\mathrm{w.dh}}_D}\mathbb{C}\,\quad\text{and} \quad\textstyle{\mathop{\mathrm{dh}}_D}\mathbb{C}\le \textstyle{\mathop{\mathrm{dh}}_{{\mathscr A}(G)}}\,\mathbb{C}\le k,
$$
where, as above, $k=\dim B$.

To obtain the first inequality note that by Part~(A) of Lemma~\ref{prochi}, the one-dimensional representation $\chi$ defined in~\eqref{ADbdet} extends to~$D$. Consider the corresponding left module $\mathbb{C}_\chi$ and its dual right module $\mathbb{C}_\chi^*$. It is shown
in \cite[p.\,643--644, Corollary~1 and the following paragraph]{Ha70}  that $\mathrm{H}^k_{\mathrm{Lie}}(\mathfrak{b},\mathbb{C}_\chi)\ne 0$. (This statement holds for any finite-dimensional Lie algebra over an arbitrary field.)
Passing from cohomology to homology,  we conclude from by the finite-dimensionality of $\mathrm{H}^{\mathrm{Lie}}_k(\mathfrak{b},\mathbb{C}_\chi)$ that
$\mathrm{H}^{\mathrm{Lie}}_k(\mathfrak{b},\mathbb{C}_\chi)\cong \mathrm{H}^k_{\mathrm{Lie}}(\mathfrak{b},\mathbb{C}_\chi^*)^*$ and hence
 $\mathrm{H}^{\mathrm{Lie}}_k(\mathfrak{b},\mathbb{C}_\chi^*)\ne 0$. Putting $M=\mathbb{C}_\chi^*$ in Theorem~\ref{toriAG}, we conclude that ${\mathop{\mathrm{Tor}}\nolimits}_k^D(\mathbb{C}_\chi^*,\mathbb{C})\ne 0$. Since $D$ is a Fr\'echet algebra, it follows from the non-triviality of ${\mathop{\mathrm{Tor}}\nolimits}_k^D(\mathbb{C}_\chi^*,\mathbb{C})$ that $k\le \mathop{\mathrm{w.dh}}_D\,\mathbb{C}$; see \cite[Proposition 2.5.4]{He00} or \cite[Proposition~4.1]{Pi08}.

Furthermore, the completion homomorphism ${\mathscr A}(G)\to D$ is a homological epimorphism by Theorem~\ref{C4mainexdi}. It follows from \cite[Proposition 3.2]{Pir_stbflat} that the application the functor $D\ptens{{\mathscr A}(G)}(-)$ to a shortest projective resolution of the ${\mathscr A}(G)$-$\mathbin{\widehat{\otimes}}$-module $\mathbb{C}$ gives a projective resolution of $D\ptens{{\mathscr A}(G)}\mathbb{C}$ in $D{\mbox{-}\!\mathop{\mathsf{mod}}}$. Since $D\ptens{{\mathscr A}(G)}\mathbb{C}\cong \mathbb{C}$ by the density of the image of ${\mathscr A}(G)$ in $D$, we conclude that the second inequality $\mathop{\mathrm{dh}}_D\,\mathbb{C}\le \mathop{\mathrm{dh}}_{{\mathscr A}(G)}\,\mathbb{C}$ holds.

To prove the third inequality we first verify that $\mathop{\mathrm{dh}}_{{\mathscr A}(G/\Lambda)}\mathbb{C}\le k$. Indeed, since $L$ is linearly complex reductive, it follows from \cite[Theorem~4.10]{Ar_aahe} that ${\mathscr A}(G/\Lambda)$ is relatively homologically trivial over ${\mathscr A}(B)$.
Proposition~\ref{beAB} implies that $0\leftarrow \mathbb{C} \leftarrow  \beta_\bullet({\mathscr A}(B),\mathbb{C})$ is a complex in ${\mathscr A}(G/\Lambda){\mbox{-}\!\mathop{\mathsf{mod}}}$. Hence by Lemma~\ref{GLaBres}, it is a projective resolution in ${\mathscr A}(B){\mbox{-}\!\mathop{\mathsf{mod}}}$. Since $\mathop{\mathrm{dh}}_{{\mathscr A}(B)}\mathbb{C}\le k$, the conditions of Proposition~\ref{BShomdi} are satisfied in the case when $A={\mathscr A}(G/\Lambda)$, $R={\mathscr A}(B)$ and $N=\mathbb{C}$. Thus, $\mathop{\mathrm{dh}}_{{\mathscr A}(G/\Lambda)}\mathbb{C}\le k$.

Finally, note that $\Lambda$ is linearly complex reductive (see \cite{Ar_lin}), and thus every projective resolution in ${\mathscr A}(G/\Lambda){\mbox{-}\!\mathop{\mathsf{mod}}}$ is a projective resolution in ${\mathscr A}(G){\mbox{-}\!\mathop{\mathsf{mod}}}$ by \cite[Proposition~5.7]{Ar_aahe}. Consequently, $\mathop{\mathrm{dh}}_{{\mathscr A}(G)}\,\mathbb{C}\le k$. Thus the third inequality is proved, which completes the argument.
\end{proof}

\end{document}